\documentstyle{amsppt}
\input amstex.tex
\loadbold

\NoRunningHeads
\magnification=1200
\pagewidth{4.98 in}
\hcorrection{.38 in}
\pageheight{7.48 in}
\vcorrection{-0.03 in}


\def\Ann {\operatorname{Ann}\nolimits}

\def\H {\operatorname{H}\nolimits}

\def\Ass {\operatorname{Ass}\nolimits}

\def\Ker {\operatorname{Ker}\nolimits}
\def\Im {\operatorname{Im}\nolimits}

\newcount\chno
\chno=0
\newcount\resno
\resno=1

\define\myref#1{\item"{[{\bf #1}]}\hfill"}

\def\proclaim #1. #2\par{\medbreak \noindent{\bf#1 \the\resno.\enspace \global\advance\resno by1} {\sl#2}\par\medbreak}

\def\xx {x_1,\ldots ,x_n}

\def\N2#1#2{\phantom{|}_{{#1}} {N}_{#2}}

\topmatter
\title Finiteness of $\bold{\bigcup_e \Ass  F^e(M)}$ and its connections to tight closure. \endtitle
\author Mordechai Katzman \endauthor
\thanks The results in this paper are part of my doctoral thesis written
at The University of Michigan, Ann Arbor,
under the supervision of Prof.~Melvin Hochster. I would like to express my deep
gratitude to Prof.~Hochster for his excellent guidance and for sharing with
me his deep mathematical insight. \endthanks
\abstract
\baselineskip=20pt
Let $R$ be a commutative Noetherian ring of prime characteristic $p$ and let
$M$ be an $R$-module.
In this paper we investigate the sets $\bigcup_e \Ass F^e(M)$ and
$\bigcup_e \Ass G^e(M)$ where $F^e()$ is the Peskine-Szpiro functor and
$G^e(M)$ is $F^e(M)$ modulo the tight closure of $0$ in $F^e(M)$.
We show that if either of these sets has finitely many maximal elements, we
can impose some restrictions on a possible counter-example
for the commutativity of localization and tight closure. The last section of this paper constructs an example which shows that $\bigcup_e \Ass F^e(M)$
 may have infinitely
many maximal elements. The question of whether $\bigcup_e \Ass G^e(M)$
may have infinitely
many maximal elements is still open.

\endabstract
\endtopmatter
\document
\baselineskip=20pt

\resno=1

\head 1. Introduction \endhead
Throughout this paper, all rings are commutative with identity and Noetherian;
$p$ will always denote a prime integer, and $q$ will be  some power $p^e$.
A local ring is defined as a Noetherian ring with a unique maximal ideal.
Let $R$ be a ring of prime characteristic $p$, let $N \subset M$ be
finitely generated $R$-modules. In \cite{HH} M.~Hochster and C.~Huneke introduced the notion of {\it the tight closure of $N$ in $M$} as follows:

Let $S$ be $R$ viewed as an $R$-algebra via the iterated Frobenius endomorphism
$r\mapsto r^q$ and define the {\it Peskine-Szpiro functor} $F^e$
from $R$-modules to $S$-modules by $F^e(M)=S \otimes_R M$. Since
the category of $S$-modules {\it is} the category of $R$-modules,
we may view $F^e$ as a functor from the category of
$R$-modules to itself.

Thus $R$-module structure on $F^e(M)$ is such that
$r^{\prime} (r \otimes m)= (r r^{\prime}) \otimes m$ and we also have
$r^{\prime} \otimes (r m)= (r^{\prime} r^q) \otimes m$.
If $I \subset R$ is an ideal then $F^e(R/I)=R/I^{[q]}$, and
generally if we apply $F^e$ to a map
$R^a \rightarrow R^b$ given by a matrix $(c_{ij})$
by identifying $F^e(R^a)\cong R^a$ and $F^e(R^b)\cong R^b$
(this identification is not canonical, it depends on a choice of
generators for the free modules),
we obtain a
map $F^e(c_{ij}): R^a \rightarrow R^b$ given by the matrix $(c_{ij}^q)$.
There is a natural map $M \rightarrow F^e(M)$ given by $m \mapsto m
\otimes 1$, and we denote the image of $m$ under this map as $m^q$.

If $N \subset M$ are $R$-modules, we have an exact sequence
$$F^e(N) \rightarrow F^e(M) \rightarrow F^e(M/N) \rightarrow 0$$
and we write $N_M^{[q]}$ for
$$\Ker ( F^e(M) \rightarrow F^e(M/N) ) \cong
\Im( F^e(N) \rightarrow F^e(M) )$$

Let $R^0$ be the set of all elements in $R$ not in any minimal prime of
$R$.
Let $N \subset M$ be $R$-modules. The {\it tight closure of N in M},
$N^*_M$, is
defined as the set of all elements $m \in M$ such that
$c m^q \in N^{[q]}_M$ for some $c \in R^0$ and all large $q$. If
$N^*_M=N$ we say that {\it $N$ is tightly closed in $M$}.

We also define $G^e(M)=F^e(M)/0^*_{F^e(M)}$. Notice that
$G^e(M/N) \cong F^e(M)/N^{[q]*}_M$.

We refer the reader to \cite{HH} for a description of the basic properties of tight closure.

\head 2. Commutativity of localization with tight closure
and the set $\bigcup_e \Ass F^e(M)$ \endhead

With the notation above, let $S\subset R$ be a multiplicative system.

We always have $S^{-1}(N^*_M) \subset (S^{-1}N)^*_{S^{-1}M}$
and we would like to know whether $S^{-1}(N^*_M)=(S^{-1}N)^*_{S^{-1}M}$.
This question still remains open in this generality. However, in a
special case, to be discussed below, an affirmative answer has been
found.

\proclaim Definition.
Let $R$ be a ring of prime characteristic $p$, and
let
$$G_\bullet  = 0 \rightarrow G_n \buildrel d_{n} \over \rightarrow
\ldots \buildrel d_{1} \over \rightarrow G_0
\rightarrow 0$$
be a complex of finitely generated projective
$R$-modules.
\item{1)} The complex $G_\bullet $ is said to have {\it phantom homology}
at the $i$th spot if $\Im d_{i+1}$ is in the tight closure of
$\Ker d_i$ in $G_i$.
\item{2)}
The complex $G_\bullet $
is said to be {\it stably phantom acyclic} if the
complex $F^e(G_\bullet )$ has phantom homology for all $i \ge 1$ and all
$e > 0$.
\item{3)} An $R$-module $M$ is said to have {\it finite phantom
projective dimension} if there exists a finite stably phantom acyclic complex
of projective $R$-modules whose zeroth homology is isomorphic to $M$.

In \cite{AHH} it is shown that if $N \subset M$ are $R$-modules and
$M/N$ has finite phantom
projective dimension then
$S^{-1}(N^*_M)=(S^{-1}N)^*_{S^{-1}M}$ for any multiplicative
system $S \subset R$, and the proof of this statement uses the fact that
under the hypotheses above $\bigcup_e \Ass F^e(M/N)$ is finite. This,
together with the following two theorems below, give a motivation for studying
the problem of whether in general $ \bigcup_e \Ass F^e (M/N) $
(or $ \bigcup_e \Ass G^e (M/N) $) is finite
or has finitely many maximal elements.

\proclaim Theorem \cite{AHH}.
Let $R$ be a ring of prime characteristic $p$, let $N \subset M$ be
finitely generated $R$-modules and let $S \subset R$ be multiplicative
system.
\item{a)} Every element of $(S^{-1}R)^0$ is a product of a unit in
$(S^{-1}R)^0$ and an element in the image of $R^0$.
\item{b)} Let $u \in M$ and $s \in S$. Then $u/w \in
(S^{-1}N)^*_{S^{-1}M}$ if and only if there exists a $d \in R^0$ such
that $s_e d u^q \in N^{[q]}_M$ for some $s_e \in S$ for all large $q$.
If $R$ has a locally stable
weak test element (resp. a completely stable weak test
element) then $d$ may be chosen to be a locally stable weak test element (resp.
a completely stable weak test element).
\item{c)} If $S$ is disjoint from $\bigcup_e \Ass F^e(M/N)$ then
 $S^{-1}(N^*_M)=(S^{-1}N)^*_{S^{-1}M}$.
\item{d)} If $S$ is disjoint from $\bigcup_e \Ass G^e(M/N)$ then
 $S^{-1}(N^*_M)=(S^{-1}N)^*_{S^{-1}M}$.

Further motivation for the study of  $\bigcup_e \Ass F^e(M)$ is given by
theorem 5 whose proof relies on the following two lemmas

\proclaim Lemma.
Let $R$ be a semi-local ring or prime characteristic $p$ with a
completely stable $q^{\prime}$-weak test element $c$, and let $J$ be
its Jacobson radical. Denote by $\widehat{\ }$ the completion with respect
to $J$.
If tight closure fails to commute with localization at
a multiplicative system $W\subset R$ for some pair
of $R$-modules $N\subset M$, then it fails to commute for a
pair of $\hat{R}$-modules.
\demo{Proof}
Pick a $u \in M$ with $u\notin N^*_M$ while
$u\in (W^{-1}N)^*_{W^{-1}M}$.
For all $q>q^{\prime}$, $c u^q \notin N^{[q]}_M$, and since $\widehat{R}$
is faithfully flat, tensoring with $\widehat{R}$ we obtain
$$c \widehat{u}^q \notin \widehat{N}^{[q]}_{\widehat{M}}$$

On the other hand, since $u\in (W^{-1}N)^*_{W^{-1}M}$, by Lemma 2b
we have $w_e c u^q \in N^{[q]}_M$ for some $w_e \in W$ and for all
$q>q^{\prime}$. Tensoring the exact sequence
$$0 @>>> \Ann_{M/N^{[q]}_M} w_e @>>> M/N^{[q]}_M @>w_e>>
M/N^{[q]}_M$$
with the faithfully flat extension $\widehat{R}$ we obtain
$({\Ann_{M/N^{[q]}_M}} w_e)\sphat =
\Ann_{\widehat{M}/ \widehat{N}^{[q]}_{\widehat{M}}} w_e$
hence $c \widehat{u}^q \in \Ann_{\widehat{M} /\widehat{N}^{[q]}_{\widehat{M}}} w_e$
and  by Lemma 2b, $u\in (W^{-1}\widehat{N})^*_{W^{-1}\widehat{M}}$.
\qed

\proclaim Lemma.
Let $R=R_1 \times \dots \times R_n$ be a ring of prime characteristic
$p$, and let $N\subset M$ be $R$-modules.
Let $N_i \subset M_i (1 \le i \le n)$ be $R$-modules, and let
$M=M_1 \times \dots \times M_n$,
$N=N_1 \times \dots \times N_n$ be the corresponding decompositions
of $N$ and $M$.
Then
$$N^*_M = {N_1^*}_{M_1} \times
\dots \times {N_n^*}_{M_n}$$
\demo{Proof}
Notice that $R^0=R_1^0 \times \dots \times R_n^0$,
$F^e_R(M) \cong F^e_{R_1}(M_1) \times \dots \times F^e_{R_n}(M_n) $
and $N^{[q]}_M= {N_1^{[q]}}_{M_1} \times \dots \times {N_n^{[q]}}_{M_n}$
where each ${N_i^{[q]}}_{M_i}$ is computed over $R_i$.
Now,
$u=(u_1,\ldots,u_n) \in N^*_M \Leftrightarrow$
there exists a $c=(c_1,\ldots,c_n)\in R^0$ such that $cu^q \in N^{[q]}$
$\Leftrightarrow c_i u_i^q \in {N_i^{[q]}}_{M_i}$ for all
$1 \le i \le n \Leftrightarrow u_i \in {N_i^*}_{M_i}$ for all
$1 \le i \le n$.
\qed

\proclaim Theorem.
Let $R$ have a $q^\prime$-weak test element, and let $N\subset M$ be
$R$-modules.
Assume that either $S=\bigcup_e \Ass F^e(M/N)$
or $S^\prime=\bigcup_e \Ass G^e(M/N)$ has finitely many maximal
elements.
If localization does not commute
with tight closure for the pair $N \subset M$, then we can find a
counter-example when $R$ is complete local
and we are localizing at a prime ideal
$P \subset R$ with $\dim(R/P)=1$.
\demo{Proof}
By Lemma 3.5a in \cite{AHH}, we may assume that we have a
counterexample in which localization at a prime
ideal $P$ fails to commute with tight closure.
Let $W$ be the complement
of $P \cup \left(\bigcup S\right)$ (respectively
$P \cup \left(\bigcup S^\prime\right)$
in $R$. By the previous theorem,
we may localize $R$ at $W$ without affecting
any relevant issues, hence we may assume that $R$ is semi-local.

Let $J$ be the Jacobson radical of $R$, and let $\widehat{}$ denote the completion
at $J$.
By Lemma 3 we can find a counter-example
over $\widehat{R}$, hence we may substitute $R$ with $\widehat{R}$,
and  we may assume that $R$ is a product of complete
local rings $R=(R_1,m_1) \times \dots \times (R_n,m_n)$
and we also get a decomposition
$M=M_1 \times \dots \times M_n$
$N=N_1 \times \dots \times N_n$ where $N_i \subset M_i$ are $R_i$ modules.

By Lemma 4 we have
$$N^*_M= {N_1^*}_{M_1} \times \dots \times {N_n^*}_{M_n}$$
and one of the pairs, say $N_1 \subset M_1$ must give a counter-example over a
local ring $(R_1,m_1)$, hence we may assume $R$ is local.

If $P=m$, we obviously cannot have a counter-example, so assume
that $P$ is not maximal, and let
$$P=P_0 \subset P_1 \subset \dots \subset P_l=m$$
be a saturated chain of primes. Let $0 \le i < l$ be the maximal
number such that localization at $P_i$ does not commute with tight closure
for the pair $N \subset M$ and replace $R,M,N$ and $P$ with
$R_{P_{i+1}}, M_{P_{i+1}}, N_{P_{i+1}}$ and $P_i$.
\qed

\proclaim Theorem.
Assume that for any local ring $(R,m)$ of prime characteristic $p$ and every
finitely generated $R$-module $\overline M$ the
set $ \bigcup_e \Ass G^e (\overline M) $ has finitely many maximal elements.
If, in addition, for every $R$-module $\overline M$ there exists a positive
integer
$B>0$ such that $m^{qB}$ kills $\H_m^0(F^e(\overline M))$
(or $\H_m^0(G^e(\overline M))$)
then tight closure commutes with localization.

\demo{Proof}
Pick a counter-example consisting of a local ring $R$, $R$-modules
$N \subset M$ and a multiplicative system $S \subset R$. By the previous
theorem we may assume that $R$ is a complete local ring, and we
are localizing at a
prime $P \subset R$ with $\dim R/P=1$. We may also assume that we have chosen
our counter-example with $\dim R$ minimal.
Since tight closure can be computed modulo the minimal primes or $R$, we
may further assume that $R$ is a domain, and hence module finite and
torsion free over a regular ring, and we may assume $R$ has weak test
elements (see section 6 in \cite{HH}.)

We may replace the pair $N\subset M$ with the pair
$0 \subset \overline M=M/N^*_M$,
hence we may assume that $N=0$ and $0$ is tightly
closed in $M$.

Pick some $u\in M$ with $u\in 0^*_{M_P}$ while $u/1 \ne 0$ in $M_P$. For all
$f\in m-P$ we have $u/1 \in 0^*_{M_f}$, otherwise we get a counter-example
over a ring of smaller dimension. Hence for all $q\gg 0$ there exists a
positive integer $N(q)$ such that $f^{N(q)} u^q=0$
in $F^e(M)$,
and since the ideal generated by all elements in $m-P$
is $m$, we have that
$ u^q\in \H^0_m(F^e(M))$ (and hence also $ u^q\in \H^0_m(G^e(M))$.)
Pick a test element $c \in R^0$.

If $m^{qB}$ kills $\H_m^0(F^e(\overline M))$, then for all
$f\in m-P$ we have
$f^{B q} c u^q=0$ and $f^B u \in 0^*_M$ but as $0^*_M=0$ we have
$f^B u =0$ in $M$, contradicting the choice of $u/1 \ne 0$ in $M_P$.

If $m^{qB}$ kills $\H_m^0(G^e(\overline M))$, then
$f^{Bq} c u^q\in 0^*_{F^e(M)}$ for all $q\gg 0$, and
by Lemma 8.16 in \cite{HH} we have $f^B u \in 0^*_M=0$ arriving again
at a contradiction.\qed
\enddemo

\head 3. The set {$\bigcup_e \Ass(F^e(M))$} over an hypersurface \endhead

In the rest of this section we will study the set $\bigcup_e \Ass(F^e(M))$,
over an hypersurface $R=A[x_1,\dots,x_n]/F$ were $A$ is a domain and with $M=R/(x_1,\dots,x_n)$. It has been shown that in some interesting cases
the set $\bigcup_e \Ass(F^e(M))$ is finite \cite{Kat}, but there is a surprisingly simple counter-example for the finiteness of this set in the general case.

\bigskip

We fix the ring $A$ to be  a domain of characteristic $p>0$,
and $q$ will always denote $p^e$ for some positive integer $e$.
For any $F\in A[\xx ]$ let $R_F={{A[\xx ]} / F}$
and let $M_F={R_F / {( \xx )R_F}}$.

\proclaim Definition.
A sequence $\{ M_n \} _n$ of $A-$modules has {\it finite torsion} if
there exists a non-zero $a \in A$ such that $(M_n)_a$ is a
torsion free $A_a$ module for all $n$.

\proclaim Lemma.
Let $A$ be a domain and let $B \supset A$ be a module finite extension domain.
Let $\{ M_i \}_i$ be a sequence of $B$ modules such that
$\{ M_i \}_i$ has finite torsion over $B$.
Then $\{ M_i \}_i$ has finite torsion over $A$.
\demo{Proof}
Pick some non zero $b\in B$ such that the modules
$(M_b)_i$ are torsion free over $B$ and
choose $a\in A$ to be a non zero multiple of $b$ in $A$.
Clearly, $(M_a)_i$ are torsion free over $A$.
\qed
\enddemo

\proclaim Lemma.
Let $A$ be a domain which is also a $k$-algebra,
let $R=A[x_1,\ldots,x_n]$ and let $I\subset R$ be an ideal generated by
elements in $k[x_1,\ldots,x_n]$.
Then any non-zero $\alpha \in A$ is a non zero divisor on $R/I$.
\demo{Proof}
Since $k[x_1,\ldots,x_n]/I$ is flat over $k$,
$R= k[x_1,\ldots,x_n]/I \otimes _k A$ is flat over $A$.
\qed

Let $A$ be a domain and let $R=A[x,y]$. If $F\in R$ is a homogeneous
polynomial, in view of Lemma 8, we may
replace $A$ with a localization at one element of a module finite extension
of $A$ to obtain
a splitting of $F$ into linear factors
$$ F(x,y)=\sum_{k=1}^s {\left( a_k x + b_k y \right )}^{r_k}$$
where $(a_k,b_k \in A$ for $1\le k \le s)$

\proclaim Lemma.
If the number of different linear factors in $F$ is
at most 3, then the modules $\{ F^e(M_F) \}_e$ have finite $A$-torsion.
\demo{Proof}
We can make a change of variables so that the different linear
factors of $F$ are among $x,y$ and $(x-y)$.

When
$F=x^{s_1}$ or $F=x^{s_1}y^{s_2}$ or $F=x^{s_1}y^{s_2}(x-y)^{s_3}$
the modules $\{ F^e(M_F) \}_e$ have no
$A-$torsion by Lemma 9.
\qed

\bigskip

In view of this lemma, the first interesting case is when $F$ is a
product of four different linear factors, and indeed our next aim
is to produce a $F$ which is a product of four linear factors
for which the modules $\{ F^e(M_F) \}_e$ {\it do not} have finite $A$
torsion.

But first we need the following lemma:

\proclaim Lemma.
Let $A$ be a domain and let $R=A[x,y]$. Then
$$ \left ( x^{q-1},y^{q-1} \right ) :_R (x-y)$$ is generated
by $y^{q-1}$ and
$\gamma=x^{q-2}+x^{q-3}y+\ldots + xy^{q-3}+ y^{q-2}$.
\demo{Proof}
Assume that $a(x-y)=bx^{q-1}+cy^{q-1}$ for some $a,b,c \in R$.
Working modulo $x-y$ we have
$\bar b x^{q-1} + \bar c x^{q-1} \equiv 0$, hence
$\bar b + \bar c \equiv 0$ and we can write
${\nobreak c=-b+d(x-y)}$ for some $d\in R$.
We can write
$ a(x-y)=bx^{q-1}+(-b+d(x-y))y^{q-1} \Rightarrow $
$ (a-dy^{q-1})(x-y)=bx^{q-1}-by^{q-1} \Rightarrow $
$ a-dy^{q-1}=b \gamma \Rightarrow$
$ a \in \left ( y^{q-1}, \gamma \right )$
\qed

\proclaim Theorem.
Let $A=k[t]$, $R=A[x,y]$. Let $F=x y (x-y) (x-ty) \in R$.
The modules $\{ F^e(M_F) \}_e$ {\it do not} have finite $A$ torsion.
\demo{Proof}
\looseness=1 We will first show that for all $q=p^e$
$$ \tau G \in \left ( x^q, y^q, xy(x-y)(x-ty) \right )$$
where $G=xy(x-y)y^{q-2}$ and $\tau=1+t+\ldots +t^{q-2}$, while\break
$G \notin \left ( x^q, y^q, xy(x-y)(x-ty) \right )$.

Let $\gamma=x^{q-2}+x^{q-3}y+\ldots + xy^{q-3}+ y^{q-2}$. We have
$\tau y^{q-2} \in (\gamma , x-ty)$ therefore
$\tau (x-y) y^{q-2} \in ((x-y)\gamma , (x-y)(x-ty))$ but by the
previous lemma, $(x-y)\gamma \in (x^{q-1},y^{q-1})$ hence
$\tau (x-y) y^{q-2} \in (x^{q-1},y^{q-1},(x-y)(x-ty))$
and
$\tau xy(x-y) y^{q-2} \in (x^q,y^q,F)$.

If $G \in \left ( x^q, y^q, xy(x-y)(x-ty) \right )$, since
$G \equiv x^2y^{q-1} \ ($mod$ (x^q,y^q))$ we can write
$x^2y^{q-1}=ax^q + by^q + cF$ for some $a,b,c \in R$ where
the $x,y$ degree of $a,b$ is $1$.
Writing $y^{q-1}(x^2-by)=ax^q+cF$ we see that $x \mid b$ and $y \mid
a$. Let $a=a^{\prime}y$ and $b=b^{\prime}x$. Notice that
now $a^{\prime},b^{\prime} \in A$. Dividing throughout by
$xy$ we get
$$y^{q-2}(x-b^{\prime}y)=a^{\prime}x^{q-1} + c(x-y)(x-ty)$$

Modulo $x-y$ this gives $y^{q-1}(1-b^{\prime})=a^{\prime}y^{q-1}$
$\Rightarrow$ $1-b^{\prime}=a^{\prime}$, while modulo $x-ty$ this gives
$y^{q-1}(t-b^{\prime})=y^{q-1} t^{q-1}a^{\prime}$
$\Rightarrow$ $t-b^{\prime}=t^{q-1}a^{\prime}$. Combining this we have
$a^{\prime}(t^{q-1}-1)=t-1$ $\Rightarrow$ $a^{\prime} \tau =1$,
which is impossible.

To finish the proof, we notice that if for some $d\in k[t]$,
$d x y (x-y) y^{q-2} \in (x^q,y^q,F)$ then
for some $a,b,c\in R$ we have
$$x \left ( dy(x-y)y^{q-2} - a x^{q-1} - cy(x-y)(x-ty)\right )=by^q$$
and since $x,y$ is a regular sequence, $x \mid b$ and $y \mid a$ and we may
write
$$d (x-y) y^{q-2} - a^\prime x^{q-1} - c(x-y)(x-ty)=b^\prime y ^{q-1}$$
where $b=x b^\prime$ and $a=y a^\prime$. Grouping together the terms
divisible by $(x-y)$ we get
$$(x-y)(dy^{q-2}-c(x-ty))\in (x^{q-1}, y^{q-1})$$
and using Lemma 11 we deduce that
$$d y^{q-2} \in (x-ty,y^{q-1},\gamma)$$
Working modulo $x-ty$ we see that
$d y^{q-2} \in (y^{q-1},\tau y^{q-2})$ hence
$\tau$ must divide $d$, and to kill all $A$-torsion we need to
invert all $\tau=\tau(q)$, and these polynomials have infinitely many
irreducible factors.
\qed

\proclaim Remark.
Notice that the counter-example above shows that the set\break
$\bigcup_e \Ass_R F_{R_F}^e(M)$ has infinitely many maximal elements.
While
$\bigcup_e \Ass F_{R_F}^e(M)$ may have infinitely
many maximal elements, the question of whether the set\break
$\bigcup_e \Ass G^e(M)$
is finite, or has finitely many maximal elements remains open.

With $R_F$ and $M_F$ as in the previous theorem, we can show that
$G^e_{R_F}(M)\cong R_F/(x,y)^q R_F$:
we can compute $0^*_{F^e(M)}$ working modulo each minimal prime of $R_F$
(see Lemma 2.10 in [AHH]). Killing the minimal primes or $R_F$ we obtain
polynomial rings, hence
$$(x,y^q)^*_{R_F/xR_F}=(x,y^q),\
(x^q,y)^*_{R_F/yR_F}=(x^q,y)$$
$$(x-y,x^q,y^q)^*_{R_F/(x-y)R_F}=
(x-y,x^q,y^q)$$
$$(x-ty,x^q,y^q)^*_{R_F/(x-ty)R_F} =(x^q,y^q,x-ty)$$ lifting these ideals
back to $R_F$
we find that $0^*_{F_{R_F}^e(M_F)}$ is the image of
$(x,y^q) \cap (x^q,y) \cap (x-y,x^q,y^q) \cap  (x^q,y^q,x-ty)$ in
$F_{R_F}^e(M_F)$.
Each monomial $x^i y^j$ is in this intersection for
all non negative integers $i,j$ with $i+j=q$, while if the image of
$H=\sum_{i+j<q} h_{ij}(t) x^i y^j$ is in $0^*_{F_{R_F}^e(M_F)}$, then
$H$ must be divisible by $x,y,(x-y)$ and $(x-ty)$, and since these are
relatively prime, $H$ must be divisible by $F$. Therefore,
$0^*_{F_{R_F}^e(M_F)}$ is the image of $(x,y)^q$ in $F_{R_F}^e(M_F)$
and $G_{R_F}^e(M_F)=R_F/(x,y)^qR_F $ and
$\Ass G_{R_F}^e(M_F)= \{ (x,y)\}$.

In fact, this argument is valid for any choice of $F$ which can be decomposed into a product of linear factors, and in view of Lemma  8, this holds
for  any homogeneous polynomial $F$.

\NoBlackBoxes

\end